\documentclass[%
 reprint,
 amsmath,amssymb,
 aps,
]{revtex4-2}
\usepackage{fontawesome5}
\usepackage{hyperref}

\usepackage{graphicx}
\usepackage{dcolumn}
\usepackage{bm}
\usepackage{xcolor}
\usepackage{amsthm}

\usepackage{slashed}

\begin{document}


\title{Independent Chiral Control in Theory-Space Models:\\
A Rank-Preserving Framework and Its Application to Neutrino Mass Generation}

\author{Aadarsh Singh}
\email{aadarshsingh@iisc.ac.in}
\affiliation{Centre for High Energy Physics, Indian Institute of Science, C V Raman Avenue, Bangalore 560 012, India}

\date{\today}

\begin{abstract}
We develop a general framework of rank-preserving, element-wise matrix
transformations for engineering fermion mass hierarchies in theory-space
constructions. We prove that preservation of massless modes requires the
transformation function to be separable, $g_f(i,j)=g^{(L)}_f(i)g^{(R)}_f(j)$, which
in turn enables independent control of left- and right-chiral zero-mode
profiles directly at the level of the theory-space mass matrix. This formalism
unifies and extends the clockwork mechanism, permits controlled deformation
of Kaluza--Klein spectra, and enhances hierarchy generation in
GIM-like fine-cancellation scenarios. As a concrete application, we show that in
theory-space models for neutrino masses, suitable transformations allow
sub-eV light neutrinos to arise from TeV-scale new physics with only
$\mathcal{O}(40)$ additional fermionic sites, while remaining consistent with
charged-lepton flavor-violation bounds. In contrast, the corresponding
untransformed models asymptote at the MeV scale and cannot access the
phenomenologically required regime without extreme field multiplicities or
hierarchical parameters.
\end{abstract}

\maketitle


\section{Introduction}
\label{sec:intro}

The origin of hierarchical scales in nature, from the electroweak scale to the
sub-eV masses of neutrinos, remains a central and unresolved puzzle in particle
physics. A particularly powerful class of Beyond the Standard Model (BSM)
theories addresses this by \emph{geometrizing} the hierarchy. In these
``theory-space'' models---including clockwork constructions and dimensional
deconstruction---small effective couplings arise not from parameters inserted by
hand, but dynamically from the localization of zero modes (massless fermionic
states) across the sites of a discretized extra dimension
\cite{giudice2017clockwork, arkani2001electroweak, kaplan2016large,
cacciapaglia2015composite, kan2024light, choudhury2022closed}. The spatial
profile of the zero mode directly determines the resulting hierarchical
suppression. In addition to geometric setups, localization may also arise from
randomness in the underlying parameters, leading to Anderson-like hierarchies
\cite{craig2018, de2020random, patel2025hierarchies}.

Despite their success, standard theory-space constructions often possess a
rigid structure. In many commonly studied models with nearest-neighbor
interactions and simple boundary conditions, the left- and right-chiral zero-mode
profiles are shaped by the same underlying Hamiltonian
\cite{grossman1999neutrino, contino2007warped, arkani2001electroweak,
csaki2003fermions, huber2001fermion, gherghetta2000bulk, chivukula2005ideal}.
As a result, the two chiralities are typically correlated, and reshaping them independently within a symmetric theory space Hamiltonian requires additional structure in conventional constructions. Independent chiral localization is a well-established feature of extra-dimensional models, here we emphasize that a similar independence can also be realized directly at the level of the theory space mass matrix through separable transformations.
Moreover, in frameworks where light masses arise from fine-cancellation (GIM-like cancellation) between
quasi-degenerate heavy states rather than from zero-mode localization, the
accessible mass scales are often restricted by the asymptotic behaviour of the
spectrum. These limitations motivate the search for a systematic method to
reshape zero-mode profiles and modify spectral features while rigorously
preserving the number of massless modes.

In this work we propose a general framework that achieves such control through a
class of index-dependent, element-wise transformations of the form
\begin{equation}
    b_{ij} = \frac{a_{ij}}{g_f(i,j)}\,.
\end{equation}
These transformations may be viewed as a generalized Hadamard (Schur) scaling
of the matrix entries \cite{Horn_Johnson_1991}, with broad relevance in areas
such as network theory, image processing, and machine learning
\cite{gonzalez2009digital, goodfellow2016deep, neudecker1995hadamard}. Our
central result establishes a necessary and sufficient condition for these
transformations to preserve the rank and nullity of the matrix, and hence the
number of zero modes. We show that massless states are preserved \emph{if and
only if} the transformation function is separable,
\begin{equation}
    g_f(i,j) = g^{(L)}_f(i)\, g^{(R)}_f(j)\,.
\end{equation}
The separability condition plays a crucial role: the factors $g^{(L)}_f(i)$ and
$g^{(R)}_f(j)$ act independently on rows and columns, allowing controlled
re-sculpting of the left- and right-chiral zero-mode profiles while keeping the
zero-mode count fixed.

This capability has significant consequences for theory-space model building.
Specific choices of $g^{(L)}_f$ and $g^{(R)}_f$ reproduce familiar constructions such as
exponential clockwork localization, warped-extra-dimensional–like profiles, or
staggered patterns, while more general choices enable smooth interpolation
between them. Crucially, the independent row/column action allows the chiral
profiles of a Dirac fermion to be modified separately---a feature that is
difficult to realize in conventional models where the zero-mode structure is
tightly constrained by the underlying Hamiltonian. In addition, these
transformations offer a systematic means of modulating the massive spectrum,
including the controlled introduction of near-degeneracies relevant for
fine-cancellation mechanisms.

We illustrate the utility of this framework through several applications.
First, we show that the standard clockwork mechanism emerges as a special case
of an appropriate separable transformation. Second, we demonstrate that
independent chiral zero-mode localization can be engineered even in symmetric
theory-space Hamiltonians. Third, we analyze how Kaluza--Klein mass towers in
deconstructed extra dimensions respond to such transformations, showing how the
massive spectrum can be reshaped without altering the number of zero modes.
Finally, we apply the framework to neutrino mass generation, where the
transformations expand the phenomenologically viable parameter space and allow
sub-eV masses to be obtained from TeV-scale physics within moderate lattice
sizes, consistent with charged-lepton-flavor-violation constraints.
It is important to emphasize the scope and intent of the present framework.
The goal is not to engineer an arbitrary fermion mass matrix with prescribed
eigenvalues and mixing matrices, which can always be achieved abstractly at the level of matrix parametrization
by writing $M = U^\dagger D V$.
Rather, we focus on identifying a restricted class of element-wise
transformations that are compatible with theory-space constructions, preserve
rank and nullity, and admit a clear physical interpretation in terms of site-local
couplings.
Within this constrained setting, the framework provides a systematic means of
reshaping localization profiles and mass spectra while maintaining the
underlying theory-space structure.

We emphasise that the framework is UV agnostic by construction; it characterizes the structural class of rank-preserving element-wise deformations available within theory space models, while the specific UV origin of any particular transformation is determined by the model under consideration. In concrete realizations such transformations can arise from physical parameters, for example, in deconstructed models, the mass matrix entries are set by site-dependent vacuum expectation values of scalar link fields, and different choices of these profiles lead to corresponding rescalings of matrix elements. Similarly, in extra-dimensional constructions, fermion masses arise from overlap integrals of left and right chiral wavefunctions, where the profiles of the two chiralities enter independently.

The remainder of this paper is organized as follows. Section
\ref{sec:framework} develops the theoretical framework for rank-preserving
transformations, including proof sketches and physical interpretation, with
complete mathematical details provided in the Appendix
\ref{app:math_framework}. Section \ref{sec:applications} presents applications
to theory-space models: generalized clockwork constructions, independent chiral
zero-mode engineering, controlled deformation of Kaluza--Klein spectra, and a
worked example in neutrino mass generation. Section \ref{sec:conclusion}
summarizes our results and discusses potential extensions.

\section{Rank-Preserving Separable Transformations of Mass Matrices}
\label{sec:framework}

Theory-space constructions encode fermion masses in matrix structures whose
null spaces correspond to massless modes in the four-dimensional effective
theory. The localization profiles of these zero modes determine the resulting
hierarchies in couplings and masses. In this section we introduce a class of
index-dependent, element-wise transformations that preserve the number of
zero modes while reshaping their profiles. This capability is the key to
engineering controlled chiral structures in BSM model building.

\subsection{Definition of the Transformation}

Consider a fermion mass matrix $A$ of dimension $m \times n$ arising from a
theory-space construction. We define an index-dependent element-wise
transformation producing a new matrix $B$ with entries
\begin{equation}
    b_{ij} = \frac{a_{ij}}{g_f(i,j)}\,,
\end{equation}
where $g_f(i,j)$ is a nonvanishing function acting on the row and column
indices. The label $f \in \mathbb{F}$ parameterizes a family of such
transformations, with $g_f : \mathbb{N} \times \mathbb{N} \to \mathbb{F}$ taking
values in a field such as $\mathbb{R}$ or $\mathbb{C}$.

This operation may be viewed as a generalized Hadamard (Schur) scaling of the
matrix entries \cite{Horn_Johnson_1991}. However, unlike a uniform rescaling,
it can in general alter the rank of $A$. Since massless fermions correspond to
zero modes of the mass matrix, our interest lies in identifying the conditions
under which the rank and nullity are preserved.

\subsection{Rank and Nullity Preservation: The Separability Condition}

The central structural property that governs rank preservation is the
separability of the transformation function. We state the main result here and
provide the detailed proof in Appendix~\ref{app:math_framework}.

\vspace{0.2cm}
\noindent\textbf{Main Result (Separability Condition).}
\textit{The transformation $b_{ij} = a_{ij}/g_f(i,j)$ preserves the rank and
nullity of $A$ if and only if $g_f(i,j)$ is separable,}
\begin{equation}
    g_f(i,j) = g^{(L)}_f(i)\, g^{(R)}_f(j)\,.
\end{equation}

\vspace{0.1cm}
\noindent
To understand the origin of this condition, consider a row $v_{i_0}$ of $A$
that is linearly dependent on the remaining rows,
\begin{equation}
    v_{i_0} = \sum_{j \neq i_0} \alpha_j v_j\,.
\end{equation}
Examining the $k$th component before and after the transformation leads to the
requirement
\begin{equation}
    \frac{g_f(j,k)}{g_f(i_0,k)} = G_f(j)\,,
\end{equation}
where $G_f(j)$ is independent of the column index $k$. This condition is
equivalent to the functional separability of $g_f$
\cite{viazminsky2008necessary}. If not separable, the transformation generally
modifies the number of linear dependencies among rows or columns, altering the
rank. Conversely, separability ensures a one-to-one correspondence between
null vectors of $A$ and $B$, guaranteeing identical rank and nullity. A complete
derivation is provided in Appendix~\ref{app:math_framework}.

\subsection{Transformation of Zero-Mode Profiles}

We now determine how the zero modes themselves transform under a separable
scaling. This result is crucial for understanding how localization profiles can
be engineered without affecting the existence of massless fermions.

\vspace{0.2cm}
\noindent\textbf{Zero-Mode Transformation.}
\textit{Let $\{v^1,\ldots,v^n\}$ be a basis for the null space of $A$. Under a
separable transformation $g_f(i,j) = g^{(L)}_f(i)\, g^{(R)}_f(j)$, the corresponding null
space of $B$ is spanned by vectors $v'^i$ whose components are}
\begin{equation}
    (v'^i)_j = (v^i)_j\, g^{(R)}_f(j)\,.
\end{equation}

This follows from substituting $a_{lj} = b_{lj}\, g^{(L)}_f(l)\, g^{(R)}_f(j)$ into the
condition $Av^i = 0$, yielding
\begin{equation}
    \sum_j b_{lj}\, (v^i_j g^{(R)}_f(j)) = 0\,,
\end{equation}
which shows that $v'^i_j = v^i_j g^{(R)}_f(j)$ spans the null space of $B$.
The column-dependent factor $g^{(R)}_f(j)$ therefore directly modifies the spatial
profile of the right-chiral zero mode, while the row-dependent factor $g^{(L)}_f(i)$
analogously controls the left-chiral profile (as relevant for Dirac fermions).
This independent row/column action enables controlled reshaping of the two
chiral sectors while preserving the zero-mode count.

\vspace{0.2cm}
\noindent\textbf{Eigenvalue Preservation.}
For completeness, we note that if $A$ is diagonalizable with eigenvalues
$\{\mu_i\}$, then $B$ shares the same eigenvalues provided the separable
function satisfies
\begin{equation}
    g^{(L)}_f(k)\, g^{(R)}_f(k) = 1 \qquad \forall\, k.
\end{equation}
This additional constraint ensures that diagonal entries transform trivially.
The proof is given in Appendix~\ref{cor:same_evalue}. Although eigenvalue
preservation is not required for our later applications, we include it as it
clarifies the broader mathematical structure of the transformation.

\vspace{0.2cm}

Having established the rank-preserving nature of separable transformations and
their action on zero-mode profiles, we now turn to explicit applications within
the context of theory-space models.

\section{Applications to Theory Space Models}
\label{sec:applications}

In this section we illustrate how rank-preserving separable transformations may
be used as a practical tool for engineering fermion mass structures in
theory-space models. Since the transformations preserve the number of zero
modes while reshaping their localization profiles, they enable controlled
manipulation of both the chiral structure and the massive spectrum. We present
several representative applications. Section~\ref{subsec:clockwork} shows that
the standard clockwork construction arises as a special case of an appropriate
separable function, while more general choices lead to qualitatively different
localization patterns. Section~\ref{subsec:chiral} demonstrates how the
independent row/column action of the transformation allows one to separate the
left- and right-chiral zero-mode profiles.  
In Section~\ref{subsec:KK} we examine the impact on Kaluza--Klein spectra in
dimensional deconstruction. Finally, in
Section~\ref{subsec:neutrino_application} we apply the framework to neutrino
mass generation as a concrete phenomenological example.

\subsection{Generalizing the Clockwork Mechanism}
\label{subsec:clockwork}

The clockwork mechanism provides a particularly elegant realization of large
hierarchies from $\mathcal{O}(1)$ parameters by generating exponentially
suppressed effective couplings through zero-mode localization in a discretized
extra dimension or theory space
\cite{giudice2017clockwork, kaplan2016large}.  
In its standard form, a nearest-neighbor structure produces a single massless
mode whose wavefunction is exponentially localized toward one end of the ``gear
chain,'' thereby inducing hierarchical interactions with fields located on
specific sites \cite{alonso2018clockwork}.  
This structure is usually introduced at the level of the Lagrangian via
site-dependent couplings or background parameters.

Our framework reproduces the clockwork construction as a special case of a
rank-preserving separable transformation and simultaneously provides a method
for exploring a broader space of localization profiles. To illustrate this,
consider the generalized clockwork Hamiltonian with site-dependent masses
$m_i$ and nearest-neighbor couplings $q_i$ \cite{ibarra2018},
\begin{equation}
    \mathcal{H}_{i,j}^{\text{cw}}
    = m_i\, \delta_{i,j} + m_i q_i\, \delta_{i+1,j}.
\end{equation}
This follows from the Lagrangian
\begin{equation}
\mathcal{L}_{\text{NP}}
    = \mathcal{L}_{\text{kin}}
    - \sum_{i=1}^n m_i \overline{L_i} R_i
    - \sum_{i=1}^n m_i q_i\, \overline{L_i} R_{i+1}
    + \text{h.c.}
\end{equation}

To isolate the role of the separable transformation, we first consider the
uniform reference case with $m_i=m$ and $q_i=1$. The corresponding Dirac mass
matrix in the basis $\{L_i,R_j\}$ is
\begin{equation}
M_{\rm CW} = m \begin{pmatrix}
1 & 1 & 0 & \cdots & 0 \\
0 & 1 & 1 & \cdots & 0 \\
\vdots & \vdots & \vdots & \ddots & \vdots \\
0 & 0 & \cdots & 1 & 1
\end{pmatrix}_{n\times(n+1)}.
\label{eq:MCW}
\end{equation}
This starting point is useful because its zero mode is delocalized. Its
components satisfy $(\xi_0)_j\propto (-1)^{j-1}$, so the magnitude of the
wavefunction is uniform across the chain, and no hierarchy is generated in this
basis.

We now apply the separable transformation
\begin{equation}
    b_{ij}=\frac{a_{ij}}{g_q(i,j)}, \qquad
    g_q(i,j)=q^{\,i-j}.
\label{eq:gq}
\end{equation}
This function can be written as
$g_q(i,j)=g_q^{(L)}(i)g_q^{(R)}(j)$ with
$g_q^{(L)}(i)=q^i$ and $g_q^{(R)}(j)=q^{-j}$. By
Corollary~\ref{cor:null_vector}, the transformed zero-mode components are
obtained by multiplying the original components by $g_q^{(R)}(j)$, giving
\begin{equation}
    \xi'_0 \propto
    \bigl\{(-q)^{-1},\,(-q)^{-2},\,\ldots,\,(-q)^{-(n+1)}\bigr\},
\end{equation}
up to an overall normalization. This is the standard exponentially localized
clockwork profile.

The transformed Hamiltonian is
\begin{equation}
    \mathcal{H}^{\text{cw}'}_{i,j}
    =
    \frac{\mathcal{H}^{\text{cw}}_{i,j}}{q^{\,i-j}}
    =
    \mathcal{H}^{\text{cw}}_{i,j}\,q^{\,j-i}.
\label{eq:Hcw_transformed}
\end{equation}
For the nonzero entries of Eq.~\eqref{eq:MCW}, this gives diagonal entries
proportional to $m$ and nearest-neighbour entries proportional to $mq$. Thus, the
separable transformation generates the canonical clockwork coupling pattern
starting from the delocalized $q_i=1$ reference matrix. Thus, the $q_i=1$ matrix provides a delocalized reference configuration, and the
separable transformation maps it to the usual clockwork form with localized
zero mode.

The framework also highlights that exponential localization is not unique.
Any separable transformation $g_q(i,j)=g^{(L)}q(i)g^{(R)}_q(j)$ preserves the zero-mode
structure while modifying its profile.  
As an illustrative example, consider
\begin{equation}
    \bar{g}_q(i,j) = q^{\,(-1)^j\, j},
\end{equation}
which depends only on the column index and therefore automatically satisfies
the rank-preserving criterion.  
The corresponding zero-mode profile becomes
\begin{equation}
    (\xi''_0)_k = (-1)^k\, q^{\,n+1 - (-1)^k k},
    \qquad k = 1,\ldots,n+1.
\end{equation}
This produces a ``staggered'' localization pattern with alternating signs and
site-dependent exponential behaviour. Both the standard and staggered profiles exhibit strong localization, as is
evident directly from the transformed zero-mode expressions, although the
functional forms governing the localization differ qualitatively. 
This demonstrates that rank-preserving separable transformations provide a
systematic method for generating qualitatively new localization patterns while
retaining the massless sector.

The method generalizes straightforwardly to non-local or random clockwork
constructions
\cite{singh2024revisitingneutrinomassesclockwork, ben2019generalized,
de2020random, patel2025hierarchies},  
where such transformations can either induce localization or enhance it.  
More broadly, separable transformations furnish a unified perspective on a wide
class of theory-space models, connecting disparate constructions through a
single mathematical framework and making explicit the structural conditions
underlying the emergence of hierarchies.

\subsection{Independent Control of Chiral Zero-Modes}
\label{subsec:chiral}
A characteristic feature of many theory-space constructions with symmetric Hamiltonians is that the left- and right-chiral zero-mode profiles are tightly correlated. For a Dirac mass matrix $M_D$, the right-chiral zero modes correspond to the null space of $M_D$, while the left-chiral zero modes are determined by the null space of $M_D^{\dagger}$. In such setups these two profiles coincide up to complex conjugation, leaving little freedom to engineer chiral asymmetries directly at the level of the theory-space mass matrix. We note that independent chiral localization is a standard feature of extra-dimensional constructions—arising, for example, from orbifold boundary conditions, bulk mass terms, or split-fermion arrangements as mentioned in introduction. The purpose of this section is to show that a similar independence of the two chiral profiles can also be achieved within the theory-space framework through separable transformations.

The separable transformations introduced in Sec.~\ref{sec:framework} provide a
systematic method to overcome this limitation. As shown in
Corollary~\ref{cor:null_vector}, a transformation of the form
$g_f(i,j) = g^{(L)}_f(i)\, g^{(R)}_f(j)$ rescales the right-chiral zero modes by the
column function $g^{(R)}_f(j)$, while the left-chiral zero modes—being null
vectors of $M_D^\dagger$—are instead rescaled by the row function $g^{(L)}_f(i)$.
Thus, the two chiralities respond to different functional components of the
same transformation. Provided $g^{(L)}_f$ and $g^{(R)}_f$ are nonvanishing, they may be
chosen independently, enabling controlled and finite reshaping of both chiral
sectors while preserving the number of zero modes.

To illustrate this capability, consider the symmetric tridiagonal Hamiltonian
\cite{tropper2021}
\begin{equation}
    \mathcal{H}_{i,j} = m\,
    (\delta_{i,j} + \delta_{i+1,j} + \delta_{i-1,j}).
    \label{Hamiltonian}
\end{equation}
For lattice sizes satisfying $n \equiv 2 \pmod{3}$, this matrix contains a
single zero mode whose components follow the repeating pattern
\[
\xi_0^k =
\begin{cases}
-1, & k \bmod 3 = 1, \\
\;\;\,1, & k \bmod 3 = 2, \\
\;\;\,0, & k \bmod 3 = 0.
\end{cases}
\]
Since the Hamiltonian is symmetric, the left- and right-chiral zero modes share the
same delocalized structure, $\vec{\xi}_{0,L}=\vec{\xi}_{0,R}\propto\{\xi^1_0,\ldots,\xi^n_0\}$.
This represents a typical situation
in which symmetric theory-space models offer no freedom to reshape the two
chiral sectors independently.

We now apply the separable transformation
\begin{equation}
    g_a(i,j) = \sin(2a\, i)\, a^{\,(-1)^j j},
\end{equation}
yielding the transformed Hamiltonian
\begin{equation}
    \mathcal{H}'_{i,j}
    = \mathcal{H}_{i,j}\, \sin(2a\, i)\, a^{\,(-1)^j j}.
    \label{new_ham}
\end{equation}
The row factor $\sin(2a\, i)$ modifies the left-chiral zero-mode profile,
while the column factor $a^{\,(-1)^j j}$ independently modifies the
right-chiral profile. The two chiralities now exhibit entirely distinct localization patterns:
The right-chiral mode acquires a staggered, exponentially varying profile
set by the column-dependent transformation, while the left-chiral mode
exhibits an oscillatory behaviour governed by the row-dependent sinusoidal
factor.

The explicit form of the transformed right-chiral zero mode follows directly
from Corollary~\ref{cor:null_vector},
\begin{equation}
    \xi'^k_{0,R} = \xi_0^k\, a^{\,(-1)^k k},
\end{equation}
while the left-chiral zero mode is given, up to normalization, by
\begin{equation}
    \xi'^k_{0,L} = \xi_0^k\, \sin(2a\, k).
\end{equation}
These expressions illustrate the general principle: separable transformations
act independently on rows and columns, and therefore allow the two chiral
sectors to be engineered separately. The particular functions chosen above are
for demonstration purposes only; any nonvanishing choice of $g^{(L)}_f(i)$ and
$g^{(R)}_f(j)$ produces a new pair of left- and right-chiral profiles while
preserving the existence of the massless mode.

For completeness, explicit closed-form expressions for the transformed
zero-mode components for different lattice sizes (including the cases
$n = 2 + 6(h-1)$ and $n = 2 + 3(2h - 1)$) are provided in
Appendix~\ref{app:chiral_examples}. These formulas illustrate how the
general scaling rule of Corollary~\ref{cor:null_vector} manifests in
specific realizations, and highlight the variety of chiral structures
that can be engineered within this framework without altering the number
of zero modes.

Although constructed as a simple example, this mechanism has broad
applicability. In flavor model building, composite Higgs frameworks, and
seesaw-type scenarios
\cite{cacciapaglia2015composite, ibarra2025flavour, covone2025flavour,
bunk2011top}, it is often advantageous to localize different chiralities in
different regions of theory space. Separable transformations provide a
systematic and model-independent method to achieve such structures starting
from otherwise rigid Hamiltonians, thereby significantly expanding the space of
realizable chiral configurations in theory-space models.

\subsection{Engineering the Kaluza--Klein Mass Spectrum}
\label{subsec:KK}

Beyond controlling zero modes, separable transformations provide a handle on the
entire tower of massive Kaluza--Klein (KK) states. Although the transformation
preserves the number of massless modes, it generically reshapes the spectrum of
massive states. Even in cases where the eigenvalues of the mass matrix are held
fixed by the condition of Corollary~\ref{cor:same_evalue}, the physical KK
masses—corresponding to singular values—can change. This distinction enables
controlled deformation of the massive spectrum without altering the zero-mode
content.

To set the stage, recall the correspondence between a five-dimensional fermion
theory and its deconstructed four-dimensional counterpart. A 5D fermion with
Lagrangian \cite{hill2002deconstructing}
\begin{equation}
\mathcal{L}_5(x^\mu,x^5)
    = \bar{\Psi}(i\gamma^\mu D_\mu - \gamma^5 D_5)\Psi
    - \frac{1}{4}\mathrm{Tr}(F_{MN}F^{MN}),
\end{equation}
compactified on an interval, yields a KK tower with masses $m_k \propto k$
determined by boundary conditions. Deconstruction replaces the 5D continuum
with a lattice of $N$ gauge groups and fermions $\Psi_k$ at sites with link
fields connecting nearest neighbors
\cite{arkani2001constructing, hill2001gauge, csaki2002exact}.  
The link-field vacuum expectation values generate a Dirac-type mass matrix whose
spectrum reproduces the KK tower in the large-$N$ limit.  
For uniform couplings of strength $M_f$, one obtains the well-known spectrum
\cite{hallgren2005neutrino}
\begin{equation}
    m_k = 2M_f \sin\!\left(\frac{k\pi}{2N}\right),
    \qquad k = 1,2,\ldots,N-1,
    \label{sin_mass}
\end{equation}
which reduces to an approximately linear spectrum for $k\ll N$,
\begin{equation}
    m_k \approx M_f\,\frac{k\pi}{N},
    \qquad
    \Delta m_k \equiv m_{k+1}-m_k \approx M_f\,\frac{\pi}{N}.
\end{equation}

Separable transformations modify this spectrum by effectively rescaling the
link-field VEVs and site-dependent masses. While a closed-form solution for a
fully arbitrary transformation is generally unavailable, useful intuition can be
gained from specific cases that are analytically tractable
\cite{kan2024light}. Consider a transformation that uniformly rescales the
diagonal and off-diagonal elements,
\begin{equation}
    g_f(i,i) = g_f, \qquad
    g_f(i+1,i) = \tilde{g}_f,
\end{equation}
leading to the transformed mass matrix
\begin{equation}
M' = M_f
\begin{pmatrix}
g_f   & 0   & 0   & \cdots & 0      \\
-\tilde{g}_f & g_f & 0   & \cdots & 0      \\
0     & -\tilde{g}_f & g_f & \cdots & 0    \\
\vdots& \vdots & \vdots & \ddots & \vdots \\
0     & 0     & 0     & \cdots & g_f \\
0     & 0     & 0     & \cdots & -\tilde{g}_f
\end{pmatrix}.
\label{trans_mass}
\end{equation}
This form may be interpreted as a separable transformation $g_f(i,j) = (\tilde{g}_f/g_f)^{i-j}$ combined with an overall
premultiplication by $g_f$. Although the eigenvalues are preserved when the
condition of Corollary~\ref{cor:same_evalue} holds, the physical KK masses are
set by the singular values of $M'$. Computing the mass-squared matrix,
\begin{equation}
\footnotesize
M'^{\!\dagger}M'
= |M_f|^2
\begin{pmatrix}
g_f^2+g_f'^2 & -g_f \tilde{g}_f & 0 & \cdots & 0 \\
-g_f \tilde{g}_f & g_f^2+g_f'^2 & -g_f \tilde{g}_f & \cdots & 0 \\
0 & -g_f \tilde{g}_f & g_f^2+g_f'^2 & \cdots & 0 \\
\vdots & \vdots & \vdots & \ddots & -g_f \tilde{g}_f \\
0 & 0 & 0 & -g_f \tilde{g}_f & g_f^2+g_f'^2
\end{pmatrix},
\end{equation}
one finds the exact eigenvalues
\begin{equation}
m_k'^2 = M_f^2\!\left(
g_f^2 + g_f'^2 - 2g_f \tilde{g}_f \cos\!\frac{k\pi}{N+1}
\right),
\qquad k=1,\ldots,N.
\end{equation}

Several features follow. In the low-$k$, large-$N$ limit,
\begin{align}
m'_k &\approx M_f |g_f - \tilde{g}_f|
    + M_f \frac{g_f \tilde{g}_f}{|g_f - \tilde{g}_f|}
      \frac{k^2\pi^2}{2N^2}, \\
\Delta m'_k &\approx
    M_f \frac{g_f \tilde{g}_f}{|g_f - \tilde{g}_f|}
    \frac{(2k+1)\pi^2}{2N^2}.
\end{align}
Thus, the transformation introduces:

\begin{itemize}
    \item a uniform offset $M_f|g_f-\tilde{g}_f|$ to all KK masses, and  
    \item a nonuniform deformation of the level spacings, making low-lying
    states relatively denser.
\end{itemize}

\begin{figure}[h!]
    \centering
    \includegraphics[width=0.4\textwidth]{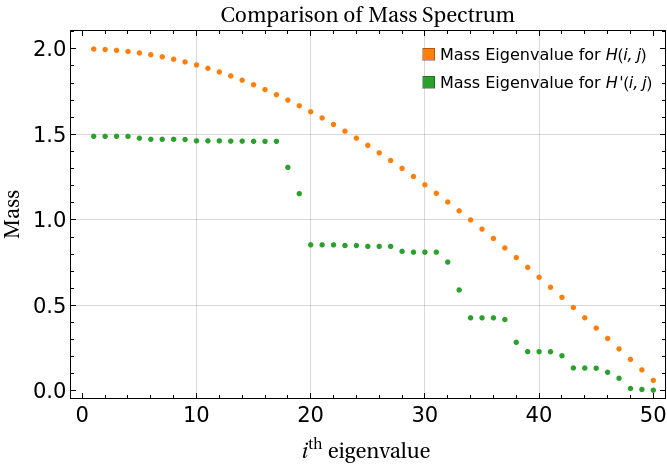}
    \caption{
    Modification of the Kaluza--Klein spectrum in a deconstructed setup.
    The orange curve corresponds to the original spectrum of
    Eq.~\eqref{sin_mass}, while the green curve shows the spectrum obtained
    from the transformed mass matrix $M'$ with
    $g_f(i,j)=(-1)^{i-j}\sin(2fi)$ (with $N=50$ and $f=1$).  
    The transformation induces regions of controlled near-degeneracy while
    preserving the number of zero modes.}
    \label{fig:Mass_Spectrum}
\end{figure}

Figure~\ref{fig:Mass_Spectrum} illustrates a more general example obtained from
the transformation $g_f(i,j)=(-1)^{i-j}\sin(2fi)$.  
Unlike uniform rescaling, the row dependence introduces site-specific
modulations that propagate through the KK spectrum.  
The resulting masses cluster into regions of near-degeneracy, a feature
frequently exploited in fine-cancellation mechanisms for neutrino physics and
other phenomenological applications
\cite{singh2024revisitingneutrinomassesclockwork}.  
We emphasise that while zero-mode localization depends exclusively on
$g^{(R)}_f(j)$ (right-chiral modes) or $g^{(L)}_f(i)$ (left-chiral modes), the massive
sector is sensitive to the combined structure of the full transformation
function.

\subsection{Application to Neutrino Mass Generation via Fine Cancellation}
\label{subsec:neutrino_application}

We now apply the transformation framework to neutrino mass generation in
theory-space models. In contrast to clockwork constructions, which generate
hierarchies via zero-mode localization, nearest-neighbor deconstruction models
with aliphatic structure typically possess no zero mode in the Dirac mass
matrix. Nevertheless, such models can yield small neutrino masses through
a fine-cancellation mechanism \cite{singh2024revisitingneutrinomassesclockwork},
where near-degenerate heavy states interfere destructively in the seesaw-like
sum that determines the effective light mass. We show that separable
transformations can enhance this cancellation by increasing the level of
near-degeneracy in the heavy spectrum, thereby enabling lighter neutrino masses
with a smaller number of sites.

The Lagrangian for the nearest-neighbor theory-space model with $n$ sites is
\cite{tropper2021}
\begin{equation}
\mathcal{L}_{\text{TS}}
    = \mathcal{L}_{\text{kin}}
    - \sum_{i,j=1}^n \bar{L}_i \mathcal{H}_{ij} R_j ,
\end{equation}
where $\mathcal{H}_{ij}$ is the symmetric tridiagonal matrix introduced in
Eq.~\eqref{Hamiltonian}. To couple the Standard Model neutrinos to the heavy fermions, we introduce
boundary interactions,
\begin{equation}
\mathcal{L}_{\text{int}}
    =
    \frac{y_1 v}{\sqrt{2}}\, \bar{\nu}_L R_n
    +
    \frac{y_2 v}{\sqrt{2}}\, \bar{\nu}_R L_1
    + \text{h.c.},
    \label{interaction-lagrangian}
\end{equation}
where $v \simeq 246~\mathrm{GeV}$ and $\langle H\rangle = v/\sqrt{2}$. 
Integrating out the heavy states at the tree level gives the light-neutrino mass
\begin{equation}
m_\nu
=
\frac{y_1 y_2 v^2}{2}
\sum_{i=1}^{n}
\frac{(u_L)^i_1 (u_R)^i_n}{m_i}\, .
\label{eq:effective_neutrino_mass}
\end{equation}
where $m_i$ are the heavy singular values, and $(u_{L/R})^i$ denote the
projections of the boundary fields onto the $i$th mass eigenstate. When the heavy spectrum becomes nearly degenerate, the sum
contains significant cancellations due to alternating phases and orthogonality,
yielding a naturally suppressed effective mass even when all underlying
parameters are $\mathcal{O}(1)$. The achievable suppression is governed by the
degree of near-degeneracy in the KK-like spectrum.

For the untransformed tridiagonal Hamiltonian, increasing degeneracy requires
a large number of sites. As shown in \cite{singh2024revisitingneutrinomassesclockwork},
masses at the MeV scale can be obtained from TeV-scale parameters with
$N \sim \mathcal{O}(10^3)$, but further suppression becomes increasingly less efficient.

A separable transformation can modify this behaviour.
As a concrete example, we consider the separable transformation
\begin{equation}
    g_f(i,j) = \sin(2f)^{\,i-j},
\end{equation}
which can be written in the separable form
$g_f(i,j)=g^{(L)}_f(i)\,g^{(R)}_f(j)$ with
$g^{(L)}_f(i)=\sin(2f)^i$ and $g^{(R)}_f(j)=\sin(2f)^{-j}$ (for $\sin(2f)\neq 0$).
This yields the transformed Hamiltonian
\begin{equation}
\mathcal{H}'_{ij}
    = \sin(2f)^{\,i-j}\, \mathcal{H}_{ij}.
\end{equation}
By construction, the transformation preserves the rank and nullity structure of the mass matrix, while
inducing a nontrivial modulation of the heavy spectrum. For definiteness, we
fix $f=1.3$ in the numerical examples shown below; this choice is representative,
and qualitatively similar behaviour is obtained for generic values of $f$
satisfying $|\sin(2f)|\neq 0,1$.
This modulation alters level spacings and can
significantly strengthen near-degeneracy across portions of the spectrum.

\begin{figure}[h!]
    \centering
    \includegraphics[width=0.4\textwidth]{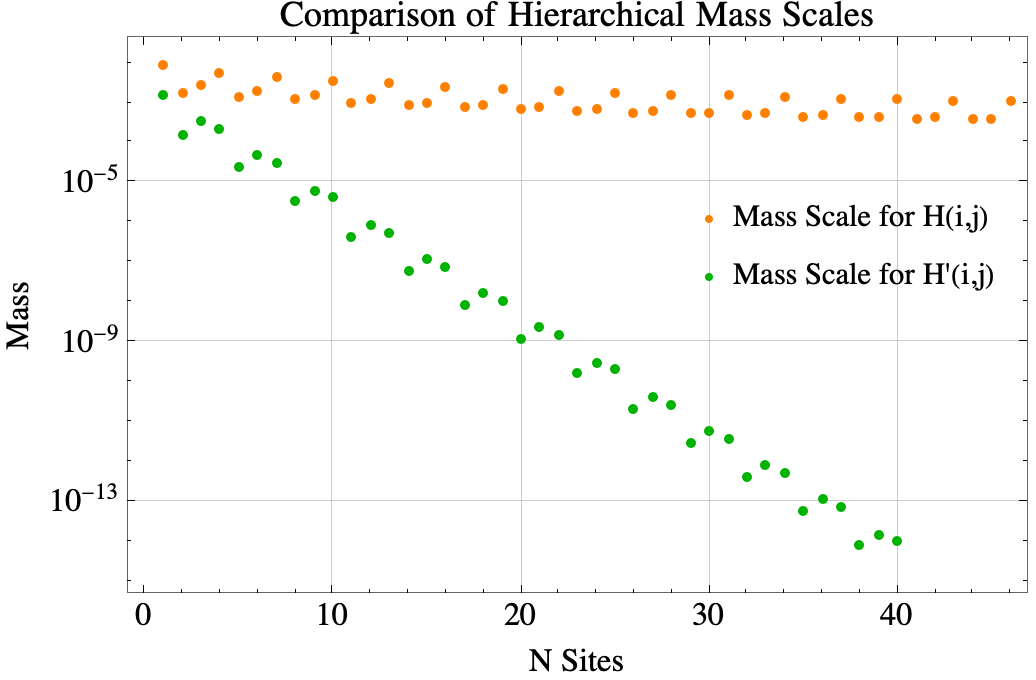}
\caption{
Achievable light-neutrino mass scale as a function of the number of theory-space
sites $N$.
Orange points correspond to the untransformed Hamiltonian $\mathcal{H}_{ij}$,
which reaches MeV-scale suppression only for $N \sim 10^{3}$ and exhibits
diminishing improvement at larger $N$.
Green points show the transformed Hamiltonian
$\mathcal{H}'_{ij} = \sin(2f)^{\,i-j}\,\mathcal{H}_{ij}$ with $f = 1.3$.
In the transformed model, the phenomenologically relevant sub-eV regime is
reached for $N \sim 30$--$40$.
}
    \label{fig:Hierarchy}
\end{figure}

Figure~\ref{fig:Hierarchy} compares the scales reached in the untransformed
and transformed cases as a function of the number of sites $N$. The
untransformed model reproduces earlier results: MeV-scale masses can be reached
with $N\sim 10^3$ for TeV new physics, but additional suppression becomes progressively inefficient.
In contrast, the transformed Hamiltonian achieves comparable suppression for
$N\sim 13$, and continues to produce substantially smaller effective masses as
$N$ increases. For $N \approx 30$--$40$, the cancellation yields values near
$0.1$~eV, consistent with the observed neutrino mass scale. The oscillatory
pattern in Fig.~\ref{fig:Hierarchy} results from the periodic structure of the
sine function.

\paragraph{Sensitivity to perturbations.}
A common concern in mechanisms that generate light masses through cancellations
among heavy modes is the sensitivity of the result to perturbations of the
underlying mass matrix. Since exact degeneracies would typically require symmetry
protection, it is important to verify that the enhanced cancellation enabled by
the transformed Hamiltonian does not hinge on finely tuned matrix entries.

To assess this, we perturb the transformed Hamiltonian according to
\begin{equation}
    \mathcal{H}' \;\rightarrow\;
    \mathcal{H}' + \delta\mathcal{H}, \qquad
    (\delta\mathcal{H})_{ij} = \epsilon_{ij}\, \mathcal{H}'_{ij},
\end{equation}
where the $\epsilon_{ij}$ are independent random variables uniformly sampled in
the range $[-\delta,\delta]$, with $\delta$ parametrizing the magnitude of the
perturbation. This corresponds to independent multiplicative perturbations of each nonzero matrix element, without enforcing any structural correlations beyond those already present in $\mathcal{H}'$. For each value of $\delta$, the effective neutrino mass $m_\nu$ is
computed via exact numerical diagonalization of the perturbed matrix, and the
resulting distribution is characterized by its median and inter-quantile range
over many realizations.

We find that, under percent-level perturbations of the transformed mass matrix,
the resulting light-neutrino mass can exhibit $\mathcal{O}(1)$ fractional
variations in regions where the cancellation is exceptionally strong and the
median approaches zero. This behaviour is expected for a mechanism based on
near-degeneracy and destructive interference among heavy states. Crucially,
however, across the entire ensemble the generated masses remain concentrated in
the same parametric regime: for $N \sim 30$--$40$, the distribution consistently
lies in the phenomenologically relevant sub-eV range, up to $\mathcal{O}(1)$
factors.

This demonstrates that the transformation robustly accesses the hierarchical
regime without requiring delicate tuning of individual matrix elements. While
the precise numerical value of $m_\nu$ is sensitive to the strength of the
cancellation, the order-of-magnitude suppression relative to the TeV-scale input
is stable under moderate perturbations. The hierarchy therefore reflects a
collective spectral property of the transformed theory-space construction rather
than an artifact of finely adjusted couplings.

These results demonstrate that separable transformations provide a controlled
method for adjusting the level spacing of the heavy spectrum without altering
the number of zero modes. In phenomenological applications such as neutrino
mass generation, where cancellations among heavy states play a central role,
this added freedom enlarges the accessible parameter space and allows realistic
neutrino masses to be obtained with moderate lattice sizes and without
introducing extremely small input parameters.

\subsection{Phenomenological Constraints}
\label{subsec:constraints}

The boundary interactions in Eq.~\eqref{interaction-lagrangian} that generate a
light neutrino mass also induce mixings between the Standard Model leptons and
the heavy theory-space states. These mixings lead to several phenomenological
signatures. Direct production of heavy neutral leptons at hadron colliders or
future lepton colliders can manifest as missing energy or displaced vertices,
depending on the mass spectrum and decay widths
\cite{abada2022collider, kwok2024searching, marcano2018heavy}. In addition,
mixing with the heavy states modifies electroweak precision observables through
loop corrections to the oblique parameters \cite{abdullahi2023present}.  
Although the fine-cancellation mechanism does not intrinsically require the new
states to lie at the TeV scale, taking them significantly below this scale with
$\mathcal{O}(1)$ couplings is strongly constrained by existing collider data.
We therefore consider masses in the multi-TeV range as a representative and
conservative regime.

A more stringent probe arises from charged-lepton flavor violation (CLFV),
especially the decay $\mu\to e\gamma$, which provides one of the strongest
constraints on neutrino mixing with heavy fermionic states
\cite{hong2019clockwork}.  
In the presence of heavy neutral leptons, the charged current becomes
\begin{equation}
\mathcal{L}_{\mathrm{CC}}
    = g W_\mu^{-} J_W^{\mu+} + \text{h.c.},
\end{equation}
with
\begin{equation}
J_W^{\mu+}
    = \sum_{\alpha=e,\mu,\tau}
      \sum_{j=1}^{n+1}
      \frac{(U_{L\alpha})^{j}}{\sqrt{2}}\,
      \bar{e}_\alpha \gamma^\mu P_L \nu_{j,\alpha},
\end{equation}
where $(U_{L\alpha})^{j}$ denotes the mixing of the SM flavor~$\alpha$ with the
$j$th mass eigenstate. These mixings enter the one-loop amplitude for
$\mu\to e\gamma$, yielding \cite{PhysRevLett.45.1908, ma1981exact}
\begin{align}
\mathrm{BR}(\mu \to e\gamma)
    &= \frac{3\alpha}{8\pi}\,|\mathcal{A}|^2,\label{eq:BRmueg}\\
\mathcal{A}
    &= \sum_{\alpha=1}^3
       \sum_{j=1}^{N+1}
       V_{\mu\alpha} V_{e\alpha}^{*}
       |(U_{L\alpha})^{j}|^2\,
       F\!\left(\frac{m_{j,\alpha}^2}{m_W^2}\right),
\end{align}
where $V_{\alpha\beta}$ is the PMNS matrix and
\begin{equation}
\footnotesize
F(x)
    = \frac{1}{6(1-x)^4}
      \left[10 - 43x + 78x^2 - 49x^3 + 4x^4 + 18x^3 \ln x\right].
\end{equation}

The MEG bound $\mathrm{BR}(\mu \to e\gamma) < 4.2\times10^{-13}$
\cite{baldini2016search} implies that, for $\mathcal{O}(1)$ Yukawa couplings,
the heavy neutrino masses must typically satisfy
$m_j\gtrsim\mathcal{O}(10)\,\mathrm{TeV}$ \cite{ibarra2018}.  
Within our framework, achieving $m_j\gtrsim 10$~TeV is compatible with the
parameter choices used in Sec.~\ref{subsec:neutrino_application}. For
$N\sim 30$--$40$ sites, the transformed Hamiltonian can simultaneously satisfy:

\begin{itemize}
    \item heavy neutrino masses in the multi-TeV regime (consistent with CLFV bounds),
    \item light neutrino masses of order $0.1$~eV (consistent with oscillation data),
    \item and $\mathcal{O}(1)$ values for the boundary Yukawa couplings.
\end{itemize}

This illustrates that the transformation enlarges the region of viable
parameter space compared with untransformed constructions, which typically
require either substantially larger $N$ or a hierarchy between diagonal and
off-diagonal mass terms.

Future experiments will continue to probe this parameter space. Heavy fermions
in the multi-TeV range remain within the prospective reach of next-generation
facilities such as the FCC and other proposed multi-TeV colliders
\cite{boyarsky2023exploring, abdullahi2023present, benedikt2020future}.  
Upcoming upgrades to CLFV searches, such as MEG II, aim to improve sensitivity
to $\mathrm{BR}(\mu\to e\gamma)$ down to $\sim 10^{-14}$
\cite{MEGII:2021fah, afanaciev2025new}, which will place even tighter
constraints on scenarios with sizable mixing between light and heavy neutrino
states.

\section{Conclusion}
\label{sec:conclusion}

In this work, we have presented a general framework for constructing fermion
mass hierarchies in theory-space models through index-dependent, element-wise
transformations that preserve the number of zero modes while reshaping their
localization profiles. The central theoretical result identifies a necessary and
sufficient condition for rank and nullity preservation: the transformation
function must be separable, $g_f(i,j)=g^{(L)}_f(i)\,g^{(R)}_f(j)$. This separability
enables controlled and independent manipulation of the left- and right-chiral
zero-mode wavefunctions, a feature that is not straightforward to realize in
conventional constructions where both chiralities are governed by the same mass
matrix structure.

We demonstrated the utility of this framework through several representative
applications. First, the standard clockwork mechanism follows as a special
case of an appropriate separable transformation, and more general choices of
$g_f$ naturally generate alternative localization patterns. Second, we showed
that independent control of the two chiral sectors can be achieved starting from
a symmetric Hamiltonian, illustrating how nontrivial chiral structures can
emerge from a simple initial configuration. Third, we analyzed the impact of
such transformations on the massive Kaluza--Klein spectrum in deconstructed
setups, showing how level spacings can be systematically reshaped and regions of
near-degeneracy introduced. Finally, we applied the framework to neutrino mass
generation via fine-cancellation, finding that suitable transformations can
enhance spectral compression and enable light-neutrino masses consistent with
oscillation data using moderate lattice sizes and $\mathcal{O}(1)$ couplings. The resulting scenarios can satisfy current charged-lepton
flavor-violation limits while remaining testable at future high-intensity and
high-energy experiments.

These results establish separable rank-preserving transformations as a
versatile tool for model building in theory space. Several natural extensions
merit further investigation. The framework is not intended as a method for arbitrary mass-matrix
construction, but as a controlled tool for exploring physically motivated
deformations within theory-space models.\\
On the phenomenological side, constructing three-generation flavor models within
this framework and exploring predictive texture structures would be of particular
interest. The present framework does not by itself predict flavor mixing angles,
since such predictions require an explicit multi-generation construction and
specified boundary couplings. Nevertheless, separable transformations act on a
Dirac mass matrix as
\begin{equation}
M \to D_L^{-1} M D_R^{-1},
\label{eq:separable}
\end{equation}
where $D_L$ and $D_R$ are diagonal but generally non-unitary matrices. Such
transformations can modify the singular-vector structure of $M M^\dagger$, and
therefore can affect the unitary matrices that enter the physical mixing matrix
in a three-generation realization. A systematic study of flavor mixing within
this framework is, therefore a natural direction for future work. A related
three-generation realization of hierarchy and flavor structure was studied in
Ref.~\cite{ibarra2025flavour}.
On the theoretical side,
a systematic classification of separable transformations and a clearer
connection to continuum descriptions—such as metric engineering in warped or
modulated extra-dimensional geometries—may reveal additional mechanisms for
hierarchy generation. Taken together, the framework developed here provides a
flexible and mathematically well-defined approach for addressing long-standing
questions related to mass hierarchies, flavor structure, and naturalness in
physics beyond the Standard Model.

\acknowledgments

A.S.\ acknowledges financial support from the Council of Scientific and
Industrial Research (CSIR), Government of India, under Senior Research
Fellowship No.\ 09/0079(15487)/2022-EMR-I.  
The author also gratefully acknowledges the use of open-source scientific
software, including Python and associated libraries, for numerical analysis.
Large Language Models (LLMs) were used solely for language polishing; all scientific content and conclusions are the author’s own.
All code used in this work is publicly available at:
\href{https://github.com/AadarshSingh0?tab=repositories}{\faGithub\ GitHub\ Profile}.

\appendix

\section{Mathematical Framework: Complete Proofs}
\label{app:math_framework}

This appendix provides formal proofs of the results stated in
Sec.~\ref{sec:framework}.  
Although our primary applications lie in fermion mass generation within
theory-space models, the underlying structure of rank-preserving element-wise
transformations is more general and may prove useful in areas such as network
theory, quantum walks, and condensed-matter eigenvalue problems.  
We therefore present the results in a theorem--proof format and adopt standard
matrix-analysis notation.

\newtheorem{theorem}{Theorem}[section]
\newtheorem{corollary}{Corollary}[theorem]
\newtheorem{definition}{Definition}[section]

\subsection{Definitions and Notation}

\begin{definition}[Index-dependent element-wise transformation]
Let $A = [a_{ij}]$ be an $N\times M$ matrix over a field $\mathbb{F}$.  
An index-dependent element-wise transformation of $A$ is a matrix
$B = [b_{ij}]$ defined by
\begin{equation}
    b_{ij} = h_f(a_{ij}, i, j),
\end{equation}
for some function $h_f:\mathbb{F}\times\mathbb{N}\times\mathbb{N}\to\mathbb{F}$.
In this work we consider the special case
\begin{equation}
    b_{ij} = \frac{a_{ij}}{g_f(i,j)},
\end{equation}
with $g_f(i,j)\neq 0$ for all $i,j$.
\end{definition}

\noindent\textbf{Notation.}
For a parameter $f\in\mathbb{F}$, we denote by
$g_f:\mathbb{N}\times\mathbb{N}\to\mathbb{F}$ a family of index-dependent
functions, where $g_f(i,j)$ acts multiplicatively on entry $a_{ij}$ of $A$.
The subscript $f$ indicates that the function may vary across a family of
transformations.  
The more general map $h_f(a_{ij},i,j)$ is used only for completeness; all
results here concern the separable case $b_{ij} = a_{ij}/g_f(i,j)$.

\subsection{Theorems and Proofs}

\begin{theorem}
\label{Them:Nullity}
Let $A$ be an $N\times M$ matrix and let $B$ be defined by
$b_{ij} = a_{ij}/g_f(i,j)$, with $g_f(i,j)\neq 0$.  
Then $A$ and $B$ have the same rank (and hence the same nullity) if $g_f$
satisfies either of the following equivalent conditions for some fixed
reference indices $(k_0,l_0)$:
\begin{equation}
    \frac{g_f(i,j)}{g_f(k_0,j)} = G_f(i)
    \qquad\text{or}\qquad
    \frac{g_f(i,j)}{g_f(i,l_0)} = \tilde{G}_f(j),
\end{equation}
where $G_f$ depends only on the row index and $\tilde{G}_f$ depends only on the column
index.
\end{theorem}

\begin{proof}
Since the transformation does not change the matrix dimensions, any trivial
nullity that arises when $N<M$ is automatically preserved.  
We therefore focus on additional linear dependencies among the rows.

Assume that the $i_0$-th row of $A$ is linearly dependent,
\begin{equation}
    v_{i_0} = \sum_{j\neq i_0} \alpha_j\, v_j,
    \label{eq:row_dependence}
\end{equation}
where $v_i$ denotes the $i$-th row of $A$.  
Let $v'_i$ denote the corresponding rows of $B$.  
Taking the $k$-th entry of Eq.~\eqref{eq:row_dependence}, we have
\begin{equation}
    v_{i_0,k} = \sum_{j\neq i_0} \alpha_j\, v_{j,k}.
\end{equation}
Using $v'_{i,k}=v_{i,k}/g_f(i,k)$, this becomes
\begin{equation}
    v'_{i_0,k}
    = \frac{1}{g_f(i_0,k)}
      \sum_{j\neq i_0} \alpha_j\, v'_{j,k}\, g_f(j,k).
\end{equation}
If the function $g_f$ satisfies
\[
    \frac{g_f(j,k)}{g_f(i_0,k)} = G_f(j),
\]
with $G_f(j)$ depending only on the row index, then the expression above becomes
\begin{equation}
    v'_{i_0,k} = \sum_{j\neq i_0} \alpha_j G_f(j)\, v'_{j,k},
\end{equation}
with coefficients $\alpha'_j = \alpha_j G_f(j)$ independent of $k$.  
Thus $v'_{i_0} = \sum_{j\neq i_0} \alpha'_j v'_j$, so the same row of $B$ is
linearly dependent.

The converse follows by the same argument applied to $B$.  
Hence $A$ and $B$ contain the same number of linearly dependent rows.  
Row–column rank equality \cite{steward1981row} then implies they also have the
same number of linearly dependent columns, so their nullities agree.  
Finally, the rank–nullity theorem gives
\[
\mathrm{rank}(A) = \mathrm{rank}(B).
\]
\end{proof}

\begin{theorem}
\label{Them:Separable}
A function $g_f(x,y)$ satisfies the conditions of Theorem~\ref{Them:Nullity}  
if and only if it is separable.
\end{theorem}

\begin{proof}
Theorem 6 of Ref.~\cite{viazminsky2008necessary} states that $g_f$ is separable
if and only if
\begin{equation}
    g_f(i,j)\, g_f(x,y)
    = g_f(x,j)\, g_f(i,y).
\end{equation}
Rearranging yields
\begin{equation}
    \frac{g_f(x,y)}{g_f(i,y)} = \frac{g_f(x,j)}{g_f(i,j)}
    \qquad\text{or}\qquad
    \frac{g_f(x,y)}{g_f(x,j)} = \frac{g_f(i,y)}{g_f(i,j)},
\end{equation}
which matches exactly the structural condition appearing in
Theorem~\ref{Them:Nullity}.  
Choosing $i$ and $j$ as reference indices identifies the corresponding
row-only and column-only functions $G_f$ and $\tilde{G}_f$, completing the proof.
\end{proof}

\subsection{Corollaries}

\begin{corollary}
\label{cor:null_vector}
Let $A$ be an $N\times M$ matrix with null-space basis
$\{v^1,\dots,v^r\}$.  
If $B$ is obtained from $A$ via a separable transformation
$b_{ij}=a_{ij}/g_f(i,j)$ with $g_f(i,j)=g^{(L)}_f(i)\,g^{(R)}_f(j)$, then a
corresponding basis for the null space of $B$ is given by
\begin{equation}
    (v'^i)_j = (v^i)_j\, g^{(R)}_f(j).
\end{equation}
\end{corollary}

\begin{proof}
Let $v^i$ be a null vector of $A$, so that
\begin{equation}
    \sum_{j=1}^{M} a_{l,j} v^i_j = 0,
    \qquad \forall\, l = 1,\dots,N.
\end{equation}
Using $a_{l,j} = b_{l,j}\, g^{(L)}_f(l) g^{(R)}_f(j)$, we obtain
\begin{equation}
    \sum_{j=1}^{M} b_{l,j}\, g^{(L)}_f(l)\, g^{(R)}_f(j) v^i_j = 0.
\end{equation}
Since $g^{(L)}_f(l)\neq 0$, it may be divided out:
\begin{equation}
    \sum_{j=1}^{M} b_{l,j}\, g^{(R)}_f(j) v^i_j = 0.
\end{equation}
Define $v'^i_j = v^i_j g^{(R)}_f(j)$. Then the above condition becomes
\begin{equation}
    \sum_{j=1}^{M} b_{l,j}\, v'^i_j = 0,
\end{equation}
so $v'^i$ is a null vector of $B$. Since the transformation multiplies each
component of each basis vector by a nonvanishing factor, the dimension of the
null space is unchanged and $\{v'^i\}$ is a basis.
\end{proof}

\begin{corollary}
\label{cor:same_evalue}
Let $A$ be an $N\times N$ diagonalizable matrix with eigenvalues
$\{\mu_i\}$ and corresponding eigenvectors $\{v^i\}$.  
Under the transformation $b_{ij} = a_{ij}/g_f(i,j)$ with
$g_f(i,j) = g^{(L)}_f(i)\, g^{(R)}_f(j)$, the transformed matrix $B$ has the same
eigenvalues as $A$ and eigenvectors
\begin{equation}
    (v'^i)_j = (v^i)_j\, g^{(R)}_f(j),
\end{equation}
if and only if
\begin{equation}
    g^{(L)}_f(k) g^{(R)}_f(k) = 1
    \qquad\forall\, k = 1,\dots,N.
\end{equation}
\end{corollary}

\begin{proof}
Let $v^i$ satisfy $A v^i = \mu_i v^i$.  
Componentwise this reads
\begin{equation}
    \sum_{j=1}^{N} a_{l,j} v^i_j = \mu_i v^i_l.
\end{equation}
Substituting $a_{l,j} = b_{l,j}\, g^{(L)}_f(l) g^{(R)}_f(j)$ gives
\begin{equation}
    \sum_{j=1}^{N} b_{l,j} g^{(L)}_f(l) g^{(R)}_f(j)\, v^i_j - \mu_i v^i_l = 0.
\end{equation}
Separate the $j=l$ term:
\begin{align}
    \sum_{j\neq l} b_{l,j} g^{(L)}_f(l) g^{(R)}_f(j) v^i_j
    + \bigl(b_{l,l} g^{(L)}_f(l) g^{(R)}_f(l) - \mu_i\bigr)\, v^i_l = 0.
\end{align}
If $g^{(L)}_f(l) g^{(R)}_f(l)=1$ for all $l$, this becomes
\begin{equation}
    \sum_{j\neq l} b_{l,j} g^{(L)}_f(l) g^{(R)}_f(j) v^i_j
    + (b_{l,l} - \mu_i) v^i_l = 0.
\end{equation}
Multiplying by $g^{(R)}_f(l)$ and defining $v'^i_j = v^i_j g^{(R)}_f(j)$ yields
\begin{equation}
    \sum_{j=1}^{N} (b_{l,j} - \mu_i \delta_l^j) v'^i_j = 0,
\end{equation}
which is the eigenvalue equation for $B$ with eigenvalue $\mu_i$.

For the converse, assume $A$ and $B$ share all eigenvalues and
$v'^i_j = v^i_j g^{(R)}_f(j)$ are eigenvectors of $B$.  
Tracing the steps above backward forces the diagonal condition
$g^{(L)}_f(k) g^{(R)}_f(k)=1$ for each $k$.
\end{proof}

\begin{corollary}
\label{cor:similarity}
If $A$ is diagonalizable and the separability condition of
Corollary~\ref{cor:same_evalue} is satisfied, then the transformed matrix
$B$ is similar to $A$.
\end{corollary}

\begin{proof}
Two diagonalizable matrices are similar if and only if they share the same set
of eigenvalues \cite{zwick2012lecture34}.  
Corollary~\ref{cor:same_evalue} ensures that $A$ and $B$ are both diagonalizable
and have identical eigenvalues.  
Hence there exists an invertible matrix $P$ such that
$B = P^{-1} A P$.
\end{proof}

\subsection{Explicit Examples of Chiral Zero-Mode Transformations}
\label{app:chiral_examples}

This subsection collects explicit closed-form expressions for the chiral
zero-mode components appearing in Sec.~\ref{subsec:chiral}.  
These formulas illustrate how the general transformation rule of
Corollary~\ref{cor:null_vector},
\[
    v'^i_j = v^i_j\, g^{(R)}_f(j),
\]
manifests in concrete theory-space Hamiltonians whose untransformed
zero-mode structure exhibits a periodic pattern.

\subsubsection*{Untransformed Zero Mode}
For the symmetric tridiagonal Hamiltonian
\begin{equation}
    \mathcal{H}_{i,j}
    = m(\delta_{i,j} + \delta_{i+1,j} + \delta_{i-1,j}),
\end{equation}
and lattice sizes satisfying $n \equiv 2 \pmod{3}$, the mass matrix possesses
a single zero mode with components
\[
\xi_{0}^k =
\begin{cases}
-1, & k \bmod 3 = 1,\\[4pt]
\;\;1, & k \bmod 3 = 2,\\[4pt]
\;\;0, & k \bmod 3 = 0,
\end{cases}
\qquad k=1,\dots,n.
\]
This mod-$3$ structure arises from the linear dependence condition among rows.

\subsubsection*{Right-Chiral Zero Mode Under the Transformation}
Consider the separable transformation
\[
    g_a(i,j)=\sin(2 a\, i)\; a^{(-1)^j j},
\]
used in Sec.~\ref{subsec:chiral}.  
The right-chiral zero mode is rescaled by the column function
$g^{(R)}_a(j)=a^{(-1)^j j}$, yielding
\[
\xi'_{0,R}{}^{\,k} = \xi_0^k\, a^{(-1)^k k}.
\]

For clarity, the explicit forms for the two periodic cases of $n$ that appear
in Sec.~\ref{subsec:chiral} are quoted below.

\paragraph{Case 1: $n = 2 + 6(h-1)$, $h\in\mathbb{N}$.}
\[
\xi_{0,R}^k =
\begin{cases}
-\,a^{\,n - (-1)^{\,n-k} k}, & k\bmod 3 = 2,\\[4pt]
\;\;\,a^{\,n - (-1)^{\,n-k} k}, & k\bmod 3 = 1,\\[4pt]
\;\;0, & k\bmod 3 = 0,
\end{cases}
\qquad k=1,\dots,n.
\]

\paragraph{Case 2: $n = 2 + 3(2h-1)$, $h\in\mathbb{N}$.}
\[
\xi_{0,R}^k =
\begin{cases}
-\;a^{-\left[n - (-1)^{\,n-k} k\right]}, & k\bmod 3 = 2,\\[4pt]
\;\;\,a^{-\left[n - (-1)^{\,n-k} k\right]}, & k\bmod 3 = 1,\\[4pt]
\;\;0, & k\bmod 3 = 0,
\end{cases}
\qquad k=1,\dots,n.
\]

These patterns differ by the effective phase assignment
$(-1)^{\,n-k}$ and the exponent sign, but both follow directly from the
multiplicative scaling rule applied to the mod-$3$ untransformed zero mode.

\subsubsection*{Left-Chiral Zero Mode Under the Transformation}
Since left-chiral zero modes correspond to the null space of
$\mathcal{H}^\dagger$ and thus pick up the row function $g^{(L)}a(i)=\sin(2 a\,i)$,
we obtain
\[
\xi'_{0,L}{}^{\,k} = \xi_0^k\, \sin(2 a\, k).
\]

Explicitly, for $n = 2 + 3h$ with $h\in\mathbb{N}$,
\[
\xi_{0,L}^k =
\begin{cases}
-\dfrac{\sin(2 n a)}{\sin(2 k a)}, & k\bmod 3 = 2,\\[10pt]
\;\;\,\dfrac{\sin(2 n a)}{\sin(2 k a)}, & k\bmod 3 = 1,\\[10pt]
\;\;0, & k\bmod 3 = 0,
\end{cases}
\qquad k=1,\dots,n.
\]

As in the right-chiral case, the structure reflects the mod-$3$
periodicity of the untransformed zero mode and the multiplicative
row-dependent scaling from the transformation.


These examples illustrate concretely how separable transformations reshape
chiral zero-mode profiles while preserving the existence of the massless state,
in agreement with Corollary~\ref{cor:null_vector}.  
They also show how analytically tractable patterns emerge when the underlying
Hamiltonian exhibits lattice periodicities, providing explicit benchmarks for
the general principles discussed in the main text.

\bibliography{apssamp}

\end{document}